\documentclass[11pt]{article}
\usepackage{latexsym}
\usepackage{amsfonts}
\makeatletter

\def\@footnotetext#1{\insert\footins{%

\footnotesize

 \interlinepenalty\interfootnotelinepenalty

 \splittopskip\footnotesep

 \splitmaxdepth \dp\strutbox \floatingpenalty \@MM

 \hsize\columnwidth \@parboxrestore

 \edef\@currentlabel{\csname p@footnote\endcsname\@thefnmark}\@makefntext
 {\rule{\z@}{\footnotesep}\ignorespaces
#1\strut}}}

\def\abstract{\small\quotation{\hskip-\parindent\sc Abstract.}}
\def\classification{\@ifnextchar [{\@xfootnotenext}%
 {\begingroup\let\protect\noexpand
 \xdef\@thefnmark{}\endgroup
 \@footnotetext}}

\title {}

\begin{document}

\classification {{\it 2000 Mathematics Subject Classification:} Primary 14E09, 
 14E25; Secondary 14A10, 13B25.\\
$\ast$) Partially supported by RGC Grant Project 7126/98P.}

\begin{center}

{\bf \Large  Non-extendable isomorphisms 
 between affine  varieties}

\bigskip

{\bf Vladimir Shpilrain} 

\medskip

 and

\medskip

 {\bf Jie-Tai Yu}$^{\ast}$

\end{center}

\medskip

\begin{abstract} In this paper, we report several large classes 
of affine  varieties (over an arbitrary field  $K$ of characteristic 0) 
 with the following property: 
 each variety in these classes has an isomorphic copy 
such that the corresponding isomorphism cannot be extended 
to an automorphism of the ambient affine space $K^n$. This implies, 
in particular, that each   of these varieties has at least two 
inequivalent embeddings in $K^n$. 
 The following application of our results seems     
interesting: we show that lines in $K^2$ are distinguished      
among  irreducible algebraic retracts by the property of having a 
unique embedding in $K^2$.  
\end{abstract} 

\bigskip

\section{Introduction }

\bigskip

 Let $K[x_1,..., x_n]$ be the polynomial algebra in $n$ variables
over a field $K$ of characteristic 0. Any collection of polynomials 
$p_1,...,p_m$ from
$K[x_1,..., x_n]$ determines an algebraic variety
${\{}p_i=0, i=1,...,m{\}}$ in the affine space $K^n$.
We shall denote
this algebraic variety by $V(p_1,...,p_m)$.
 \smallskip

 We say that two algebraic varieties $V(p_1,...,p_m)$ and $V(q_1,...,q_k)$
 are {\it isomorphic} if the algebras of residue
classes $K[x_1,..., x_n]/\langle p_1,...,p_m \rangle$ and
$K[x_1,..., x_n]/\langle q_1,...,q_k \rangle$ are
isomorphic.
Here $\langle p_1,...,p_m \rangle$ denotes the ideal of
 $K[x_1,..., x_n]$ generated by $p_1,...,p_m$.
 Thus, isomorphism that we consider here is algebraic, not geometric, 
 i.e., we actually consider isomorphism of what is called  {\it affine schemes}. 

 On the other hand, we say that two algebraic varieties
 (or, rather, embeddings of the same algebraic variety in $K^n$) are
{\it equivalent} if there is an automorphism of $K^n$ that takes one
of them onto the other. Algebraically, this means there is 
an automorphism of $K[x_1,..., x_n]$ that takes the ideal 
$\langle p_1,...,p_m \rangle$ to the ideal $\langle q_1,...,q_k \rangle$. 
 \smallskip 

 Given two isomorphic  varieties $V_1$,  $V_2$ and an isomorphism 
$\varphi : V_1 \to V_2$ between them, this isomorphism may or may not be 
 extended to an automorphism of the ambient affine space $K^n$. 
(In a purely algebraic language, i.e., when talking about isomorphism 
between algebras of residue classes, we would say that the isomorphism 
may or may not be  {\it lifted} to an automorphism of $K[x_1,..., x_n]$.) 
 Clearly, if an isomorphism $\varphi : V_1 \to V_2$ cannot be 
extended to an automorphism of $K^n$, then the varieties $V_1$ and $V_2$ 
represent two inequivalent embeddings of the same algebraic variety in $K^n$. 
 \smallskip

In this paper, we concentrate on examples of varieties with non-extendable 
isomorphisms, and therefore with inequivalent embeddings in $K^n$. 
Numerous examples of that sort were published previously -- see e.g. \cite{AS}, 
\cite{Harm}, \cite{K},   \cite{ShpYu99}, \cite{ShpYuHyp}. However, all 
varieties in those examples required an individual approach, i.e., 
there was, up until now, no reasonably big ``pool" of such examples. 
 In particular, in our earlier paper \cite{ShpYuHyp}, even though we have 
reported a very simple criterion for detecting varieties with 
inequivalent embeddings, it could not be decided just by inspection whether 
a given variety satisfies this criterion; again, each example had to be 
treated individually. 

 Here we are able to give a couple of  criteria for an algebraic variety 
to have inequivalent embeddings,  that  
can be verified just by inspection, at least in the case where a 
given variety is a hypersurface, i.e., is of the form $V(p)$ for some 
polynomial $p$. These criteria, in general, appeal to a Gr\"obner basis (a 
general reference here is \cite{AL}) 
of the ideal that determines a given variety; however, in the case where 
a variety is a hypersurface $V(p)$, the criteria can be formulated without 
mentioning  Gr\"obner bases, but just by appealing to the collection of 
monomials of the polynomial $p$.  Crucial for our criteria is the following 
result of Hadas \cite{Hadas}: if $q$ is a {\it coordinate} polynomial of 
$K[x_1,..., x_n]$ (i.e., $q$ can be taken to $x_1$ by an automorphism of 
 $K[x_1,..., x_n]$), then all vertices of the Newton polytope of $q$ are 
on coordinate hyperplanes. This condition can be easily translated into a 
collection of simple linear inequalities for the exponents of variables to 
which they occur in monomials of $q$. 

 Recall that a pure lexicografic order on the set of monomials in 
the variables $x_1,..., x_n$ is induced by an order $x_{i_1} < ...< x_{i_n}$
on the set of variables, so that, for example, $x^k_{i_1} < x_{i_n}$ for 
any $k$. 

 Now we are ready to formulate our criteria. 
\medskip

\noindent {\bf Theorem 1.1.} Let $p_1=p_1(x_1,..., x_n)=x_1-f(x_1^k, x_2,...,x_n),
 ~k \ge 2$, ~where $x_1$ actually occurs in the polynomial $f$. 
Let $p_i=p_i(x_2,...,x_n), ~2  \le i \le m$. Suppose every polynomial in 
 the Gr\"obner basis, with respect to some pure lexicografic order, 
 of the ideal $\langle p_1,...,p_m \rangle$, has the highest monomial of 
the form $x_1^{k_1} \cdot ... \cdot x_n^{k_n}$, with all  $k_i > 0$. 
Then the isomorphism ~$\varphi : x_1 \to x_1^k, ~x_i \to x_i, ~i \ge 2$, 
~between $K[x_1,..., x_n]/\langle x_1-f^k(x_1, x_2,...,x_n), p_2,..., p_m 
\rangle$ and $K[x_1,..., x_n]/\langle p_1,...,p_m \rangle$, cannot 
be lifted to an automorphism of $K[x_1,..., x_n]$. 
\medskip

 We have to admit that the condition on the Gr\"obner basis 
in Theorem 1.1 can be computationally hard to verify if $m \ge 2$, and, moreover, 
a Gr\"obner basis is not very likely to have the required form in that case. 
However, for a single  polynomial $p_1$, mentioning a Gr\"obner basis 
can be altogether  avoided, and the criterion becomes simple:

\medskip

\noindent {\bf Corollary 1.2.} Let $p=p(x_1,..., x_n)=
x_1-f(x_1^k, x_2,...,x_n),
 ~k \ge 2$, ~where $x_1$ actually occurs in the polynomial $f$. 
Suppose the Newton polytope of  the polynomial $f$ has a vertex outside of 
any coordinate hyperplane. Then the isomorphism ~$\varphi : x_1 \to x_1^k, 
~x_i \to x_i, ~i \ge 2$, ~between $K[x_1,..., x_n]/\langle 
x_1-f^k(x_1, x_2,...,x_n) \rangle$
 and $K[x_1,..., x_n]/\langle p \rangle$ cannot 
be lifted to an automorphism of $K[x_1,..., x_n]$. 
\medskip

 Another large class  of algebraic varieties with inequivalent embeddings 
is given in the following 
\medskip

\noindent {\bf Theorem 1.3.} Let $p_1=p_1(x_1,..., x_n)=
x_1-f_1(x_1 \cdot f_2, x_2,...,x_n)$, ~where $x_1$
 actually occurs in the polynomial $f_1$, and the polynomial 
$f_2=f_2(x_2,..., x_n)$ is  not a constant.
 Let $p_i=p_i(x_2,...,x_n), ~2  \le i \le m$. Suppose every 
polynomial in 
 the Gr\"obner basis, with respect to some pure lexicografic order, 
 of the ideal $\langle p_1,...,p_m \rangle$, has the highest monomial of 
the form $x_1^{k_1} \cdot ... \cdot x_n^{k_n}$, with all  $k_i > 0$, and 
 suppose that every such  monomial is higher than any monomial in $x_1 \cdot f_2$. 
Then the isomorphism ~$\varphi : x_1 \to x_1 \cdot f_2, ~x_i \to x_i, ~i \ge 2$, 
~between 
$K[x_1,..., x_n]/\langle x_1-f_1(x_1, ...,x_n) \cdot f_2, p_2, ..., p_m \rangle$
 and $K[x_1,..., x_n]/\langle p_1, p_2,..., p_m 
\rangle$ ~cannot be lifted to an automorphism of $K[x_1,..., x_n]$. 
\medskip

 Again, for a single  polynomial $p_1$, the criterion becomes very simple:

\medskip

\noindent {\bf Corollary 1.4.} Let $p=p(x_1,..., x_n)=
x_1-f_1(x_1 \cdot f_2, ~x_2,...,x_n)$, ~where $x_1$
 actually occurs in the polynomial $f_1$, and the polynomial 
$f_2=f_2(x_2,..., x_n)$ 
is  not a constant. Suppose the Newton polytope of  the polynomial 
$p$ has a vertex outside of 
any coordinate hyperplane, and suppose that $p$ has higher monomials  
 than any monomial in $x_1 \cdot f_2$. 
Then the isomorphism 
~$\varphi : x_1 \to x_1 \cdot f_2, ~x_i \to x_i, ~i \ge 2$, 
~between $K[x_1,..., x_n]/\langle x_1-f_1(x_1, ...,x_n) \cdot f_2(x_2,...,x_n) \rangle$ and 
$K[x_1,..., x_n]/\langle p \rangle$ cannot 
be lifted to an automorphism of $K[x_1,..., x_n]$. 
\medskip

 This has another interesting corollary, which is a rather unexpected 
payoff of our method. Recall that a subalgebra 
$S$ of an algebra $R$ is called a 
 {\it  retract} if there is an idempotent homomorphism (a {\it retraction}, 
 or {\it projection})  $\varphi:   R \to  R$    such 
  that   $\varphi(R) = S.$ 
A characterization of retracts of a two-variable 
polynomial algebra $K[x, y]$ was given in \cite{ShpYu00}. Since every 
proper retract of $K[x, y]$ is of the form $K[p]$ for some polynomial 
$p=p(x,y)$ (see \cite{Costa}), we shall also call the curve $p=0$ a (algebraic)  
retract of $K^2$ if $K[p]$ is a retract of $K[x, y]$. 
Based on our characterization of retracts \cite{ShpYu00} 
 and on Corollary 1.4 above, we get

\medskip

\noindent {\bf Corollary 1.5.} Let $p=p(x,y)$, and let $K[p]$ be a 
retract of $K[x, y]$. The curve $p=0$ has  inequivalent embeddings in $K^2$ 
unless $p$ is either  a coordinate polynomial or can be taken 
to $xy$ by an automorphism of $K[x, y]$. 
\medskip

The meaning of this result is that it distinguishes  lines  
among irreducible  algebraic retracts of $K^2$ by means of an ``external" 
property of having a unique embedding in $K^2$. The 
fact that $p=0$ has a unique embedding in $K^2$ for a coordinate polynomial $p$,  
is a well known result of Abhyankar and Moh 
\cite{AM}  and  Suzuki \cite{Suzuki}. The  curve $xy=0$ is known to have 
a  unique embedding in ${\bf C}^2$ -- see  \cite{J2}.  

\medskip

In  Section 3, we consider embeddings 
of    varieties of codimension 2. 
  It is known that every algebraic variety in ${\bf C}^n$  
 has a unique embedding in ${\bf C}^{2n+2}$ -- 
see  \cite{K}  or \cite{Srinivas}. The situation with embeddings 
of  varieties of a high codimension  is therefore  really 
intriguing. For instance, to the best of our knowledge, there is no   example  
of  a smooth irreducible algebraic variety of dimension $n$ with 
inequivalent embeddings in ${\bf C}^{2n+1}$.  
 
\smallskip

 Here we give examples, for any $n \ge 3$, of algebraic varieties 
of codimension two having inequivalent embeddings in $K^n$. These varieties 
however are not smooth, even though each has only  one singular  point.  
We note that  Kaliman \cite{K}  gave an example of a  curve (with one singular
 point) that has inequivalent embeddings  in  ${\bf C}^3$. Our method here 
seems to be more ``generic", i.e., we, in fact, give a rather general recipe for 
constructing examples of that sort. \\

\section{Varieties with non-extendable isomorphisms  }

\bigskip

\noindent {\bf Proof of Theorem 1.1.} The fact that the mapping $\varphi$ 
is actually an isomorphism, is explained in \cite[Example 1]{ShpYuHyp}. 
We  are not going to reproduce the argument here  because this would 
require too much background material. However, to make the exposition 
here as self-contained as possible,  we  verify that $\varphi$ 
is an  {\it onto} homomorphism. The fact that it is onto is fairly 
obvious since, modulo the ideal $\langle p_1,...,p_m \rangle$, we have 
$x_1=f(x_1^k, x_2,...,x_n)$.
 To see that $\varphi$ is a homomorphism, observe that 
$\varphi(x_1-f^k(x_1, x_2,...,x_n))= x_1^k - f^k(x_1^k, x_2,...,x_n)=
(x_1-f(x_1^k, x_2,...,x_n))\cdot (...)$. 

 Assume now, by way of contradiction, that $\varphi$ can  
be lifted to an automorphism of $K[x_1,..., x_n]$. Then there must be 
a coordinate polynomial of the form $x_1^k + u(x_1,..., x_n)$, where 
the polynomial $u=u(x_1,..., x_n)$  belongs to the 
ideal $\langle p_1,...,p_m \rangle$. The highest monomial 
(with respect to a given pure lexicografic order) of the 
polynomial $u$ has to be therefore divisible by the highest monomial 
of some polynomial in the Gr\"obner basis
 (with respect to the same pure lexicografic order) of the 
ideal $\langle p_1,...,p_m \rangle$. We claim that, if this is the case, 
then the polynomial $x_1^k + u$ has a vertex outside of 
any coordinate hyperplane, and therefore cannot be coordinate. 
 Indeed, if all vertices of $x_1^k + u$ were on coordinate hyperplanes, 
that would mean that for any monomial $x_1^{k_1} \cdot ... \cdot x_n^{k_n}$
 of $u$ involving all variables, there is another monomial 
$x_1^{m_1} \cdot ... \cdot x_n^{m_n}$, involving at least one variable 
less, such that $m_i > k_i$  for at least one $i$. Then the highest monomial 
of $x_1^k + u$ with respect to (any) pure lexicografic order would 
have the same form, i.e., would be missing at least one variable. A 
monomial like that cannot be divisible by 
$x_1^{k_1} \cdot ... \cdot x_n^{k_n}$ ~with all  $k_i > 0$, hence a 
contradiction. $\Box$ 
\medskip

\noindent {\bf Proof of Theorem 1.3} is similar. Again, the fact that the 
mapping $\varphi$ 
is actually an isomorphism, is explained in \cite[Example 2]{ShpYuHyp}. 
 To see that $\varphi$ is a homomorphism, observe that 
$\varphi(x_1-f_1(x_1, ...,x_n) \cdot f_2(x_2,...,x_n)) = 
x_1 \cdot f_2 - f_1(x_1 \cdot f_2 , ...,x_n) \cdot f_2 = 
(x_1 - f_1(x_1 \cdot f_2 , ...,x_n)) \cdot f_2$. Then, $\varphi$ is obviously 
onto since, modulo the ideal 
$\langle x_1-f(x_1 \cdot f_2, x_2,...,x_n),...,p_m \rangle$, we have 
 $x_1=f(x_1 \cdot f_2, , x_2,...,x_n)$. 

 The proof of the fact that $\varphi$ cannot be lifted to an automorphism 
of $K[x_1,..., x_n]$ goes   along   the same lines as in 
the proof of Theorem 1.1, with one  difference. The additional restriction 
  in the statement of Theorem 1.3 is needed because, if some  monomial of 
$x_1 \cdot f_2$  is equal to the highest
monomial of some polynomial  in the Gr\"obner basis  
 of the ideal $\langle p_1,...,p_m \rangle$, then, due to cancellations,  
  the  polynomial  ~$\varphi(x_1)$ may be lifted to a  coordinate 
of $K[x_1,..., x_n]$. This will be important to us in the proof of 
Corollary 1.5. $\Box$   

\medskip

\noindent {\bf Proof of Corollary 1.5.} If $K[p]$ is a 
retract of $K[x, y]$, then, by   \cite[Theorem 1.1]{ShpYu00}, 
 there is 
 an automorphism $\psi$ of $ K[x, y]$ that takes the  polynomial 
$p$ to $q(x,y)=x + y \cdot f(x,y)$ ~for some polynomial $f(x,y)$. Since 
the curve $p=0$ has a unique embedding in $K^2$ 
if and only if $q=0$ does, we may as well assume that $p$ itself is 
of the form $p(x,y)=x + y \cdot f(x,y)$. Then, by Corollary 1.4 (with 
$f_1(x,y)=-f(x,y), ~f_2(x,y)=y$), the curve $p=0$ has inequivalent 
embeddings in $K^2$ unless the degree in  $x$ of 
 the polynomial $f(x,y)$ is $\le 1$. If the degree is 0, then 
 $p(x,y)=x + \sum a_i y^i$ 
is a coordinate polynomial, and   therefore has a unique embedding in $K^2$ 
by  a well known result of Abhyankar and Moh \cite{AM}  and  Suzuki \cite{Suzuki}. 

 If the degree is 1, then Corollary 1.4 is applicable unless 
$f(x,y)=x+g(y)$. In that case,  $p(x,y)=x + xy+ y\cdot g(y)$. Apply the 
automorphism $x \to x, ~y  \to y-1$. Then $p(x,y)$ becomes 
$p'(x,y)= xy+(y-1)\cdot g(y-1)$. This can be written as $y\cdot (x+h(y)) +c$, 
where $c \in K$. After applying the automorphism $x \to x -h(y), ~y  \to y$, 
this becomes $xy +c$. Now if $c \ne 0$, the curve $xy +c=0$ is well known to 
have inequivalent embeddings in $K^2$. If $c =0$, we have the curve $xy =0$ 
that has a unique embedding in ${\bf C}^2$ by a result of Jelonek \cite{J2}.
$\Box$  \\

\section{Varieties  of   codimension two}

\bigskip

\noindent {\bf Example 3.1.} Let the curve $C_1$ in, say, ${\bf C}^3$, 
 be the common zero locus of two polynomials, $p_1= p_1(x,y,z)= 
x-x^2y-yz-z+ \frac{1}{4}$ 
and $q_1= q_1(x,y,z)= y-z^2- \frac{1}{2} z + 2xy +2x- \frac{15}{16}$. 
 The gradients of $p_1$ and $q_1$ have the only common zero at the 
point $(-\frac{1}{2}, -1, - \frac{1}{4})$, and, since this point belongs to 
the curve $C_1$, it is a singular point of this curve. (This is, in fact, 
the only singular point of $C_1$.)  
 Therefore, if a curve $C_2$, which is the common zero locus of two 
polynomials $p_2$  and $q_2$, is to be 
 equivalent to $C_1$ under an automorphism of ${\bf C}^3$, then there 
should be a point where both the gradients of 
$p_2$ and $q_2$ are equal to 0. (This follows easily from the ``chain rule" 
for partial derivatives.) 
 We are now going to exhibit a curve $C_2$ which is isomorphic to $C_1$ 
but has no   points of this kind.

As in \cite{ShpYu99}, it  will be technically more 
convenient  to write   algebras of residue classes as 
``algebras with relations", i.e., for example, instead of 
${\bf C}[x_1,...,x_n]/\langle p(x_1,...,x_n) \rangle$ we shall write 
$\langle x_1,...,x_n \mid p(x_1,...,x_n)=0\rangle$. 

 Now we  get the following chain of ``elementary" isomorphisms:
\smallskip

\noindent  $ \langle x,y,z \mid x= x^2y+yz+z- \frac{1}{4}, 
~y =  z^2+ \frac{1}{2} z - 2xy -2x+ \frac{15}{16} \rangle  ~\cong \\
 \langle x,y,z,u \mid u=xy, ~x=xu+ yz+z- \frac{1}{4}, 
~y =  z^2+ \frac{1}{2} z - 2u -2x+ \frac{15}{16} \rangle ~\cong \\ 
\langle x,z,u \mid u= x z^2+ \frac{1}{2} xz - 2xu -2x^2+ \frac{15}{16}x, 
  ~x=xu+ z- \frac{1}{4} + (z^2+ \frac{1}{2} z - 2u -2x+ \frac{15}{16})z \rangle 
~\cong \\
 \langle x,y,z \mid y= x z^2+ \frac{1}{2} xz - 2xy -2x^2+ \frac{15}{16}x, 
 ~x=xy+ z- \frac{1}{4} + (z^2+ \frac{1}{2} z - 2y -2x+ \frac{15}{16})z \rangle$. 

 Thus, we have a curve $C_2$ which is isomorphic to $C_1$ and 
which is the common zero locus of $p_2= p_2(x,y,z)= 
y- x z^2- \frac{1}{2} xz + 2xy +2x^2- \frac{15}{16}x$  
and $q_2= q_2(x,y,z)= x-xy- z+ \frac{1}{4} - 
(z^2+ \frac{1}{2} z - 2y -2x+ \frac{15}{16})z$.  

 The  gradient of $p_2$    vanishes only at the point 
$(-\frac{1}{2}, \frac{5}{4}, - \frac{1}{4})$, whereas the gradient of $q_2$ 
does not vanish at this point, i.e., the gradients of $p_2$ and $q_2$ 
have no common zeros. Therefore, $C_2$ is not equivalent to $C_1$.

\medskip

\noindent {\bf Example 3.2.} Based on the previous example, we can
construct examples of algebraic varieties of codimension two with 
inequivalent embeddings in $K^n$ for any $n \ge 3$ as follows. 
 Let $n = 3+k, ~k \ge 1$, and let $p_1= p_1(x,y,z, t_1, ..., t_k)= 
x-x^2y-yz-z+ \frac{1}{4} + t_1^2 + ... + t_k^2$; ~$q_1= 
q_1(x,y,z, t_1, ..., t_k)= y-z^2- \frac{1}{2} z + 2xy +2x- \frac{15}{16}$. 
Let $V_1$ be the common zero locus of $p_1$ and $q_1$ in $K^n$. 
 Then the only singular point of $V_1$ where both the gradients of 
$p_1$ and $q_1$ are equal to 0, is the point 
$(-\frac{1}{2}, -1, - \frac{1}{4}, 0, ..., 0)$. 

 Arguing as in Example 3.1, we get a variety $V_2$ which is isomorphic to 
$V_1$ and 
which is the common zero locus of $p_2= p_2(x,y,z, t_1, ..., t_k)= 
y- x z^2- \frac{1}{2} xz + 2xy +2x^2- \frac{15}{16}x$  
and $q_2= q_2(x,y,z, t_1, ..., t_k)= x-xy- z+ \frac{1}{4} - 
(z^2+ \frac{1}{2} z - 2y -2x+ \frac{15}{16})z + t_1^2 + ... + t_k^2$. 
The gradients of $p_2$ and $q_2$ 
have no common zeros, hence $V_2$ is not equivalent  to $V_1$.  
\medskip

 We note that the choice of constant terms  in the polynomials $p_1$ 
and $q_1$ was made so that the point where both the gradients of 
$p_1$ and $q_1$ are equal to 0 would belong to our variety $V_1$. We needed 
this to be able to prove that $V_2$ is not equivalent  to $V_1$.  
 However, it seems plausible (although we do not have a proof 
at this time) that (the corresponding) $V_2$ is not going to 
be equivalent  to $V_1$ 
with most   any choice of constant terms  in $p_1$ and $q_1$, which would, 
in particular,   give an example of a smooth irreducible curve with 
inequivalent embeddings in ${\bf C}^3$. \\

\noindent 
 Department of Mathematics, The City  College  of New York, New York, 
NY 10031 
\smallskip

\noindent {\it e-mail address\/}: ~shpil@groups.sci.ccny.cuny.edu 
\smallskip

\noindent {\it http://zebra.sci.ccny.cuny.edu/web/shpil} \\

\noindent Department of Mathematics, The University of Hong Kong, 
Pokfulam Road, Hong Kong 
\smallskip

\noindent {\it e-mail address\/}: ~yujt@hkusua.hku.hk 
\smallskip

\noindent {\it http://hkumath.hku.hk/\~~jtyu}

\end{document}